\documentclass[12pt,letterpaper]{article}
\usepackage{amsfonts}
\usepackage{amssymb}
\usepackage{latexsym}  
\usepackage{pifont}  

\newcommand\DD{{\mathcal D}}

\newcommand\FF{{\mathcal F}}
\newcommand\GG{{\mathcal G}}
\newcommand\HH{{\mathcal H}}
\newcommand\NN{{\mathcal N}}
\newcommand\PP{{\mathcal P}}

\renewcommand\SS{{\mathcal S}}

\newcommand\CC{{\mathcal C}}
\newcommand\UU{{\mathcal U}}

\newcommand\RRR{{\mathbb R}}
\newcommand\ZZZ{{\mathbb Z}}
\newcommand\QQQ{{\mathbb Q}}
\newcommand\TTT{{\mathbb T}}

\newcommand\cccc{\mathfrak{c}}
\newcommand\pppp{\mathfrak{p}}

\newcommand\res{\mathord {\upharpoonright}}  

\newcommand\Hom{\mathrm{Hom}}   

\newcommand\bohr{\mathsf{b}} 

\newcommand\iv{^{-1}} 

\newcommand\mstar{{\raise -1 pt \hbox{\ding{99}}}}
\newcommand\mdag{{\raise -1 pt \hbox{\ding{63}}}}

\newcommand\up{{\fam0 \uparrow}}

\newcommand\NNb{{\NN_{\mathsf{b}}}}

\newcommand\eop{$\ \ {\vcenter
   {\hrule
   \hbox{\vrule height 9pt \kern 9pt \vrule height 9pt}
   \hrule}}$\vskip 1.0 pt}

{\end{enumerate}}
\newenvironment{itemizz}{\begin{itemize}\setlength{\itemsep}{-1mm}} %
{\end{itemize}}

\newtheorem{theorem}{Theorem}[section]
\newtheorem{definition}[theorem]{Definition}
\newtheorem{lemma}[theorem]{Lemma}
\newtheorem{corollary}[theorem]{Corollary}
\newtheorem{proposition}[theorem]{Proposition}
\newtheorem{example}[theorem]{Example}
\newtheorem{remark}[theorem]{Remark}

\newenvironment{proof}{{\bf Proof.}}{\eop\medskip}
\newenvironment{proofof}[1]{\medskip \textbf{Proof of #1.}}{\eop\medskip}

\begin{document}

\title{
Limits in Function Spaces and Compact Groups}

\author{Joan E. Hart\footnote{University of Wisconsin, Oshkosh,
WI 54901, U.S.A.,
\ \ hartj@uwosh.edu}
\  and
Kenneth Kunen\footnote{University of Wisconsin,  Madison, WI  53706, U.S.A.,
\ \ kunen@math.wisc.edu}
\thanks{Both authors partially supported by NSF Grant DMS-0097881.}
}

\maketitle

\begin{abstract}
If $B$ is an infinite subset of $\omega$ and $X$ is a topological group,
let $\CC^X_B$ be the set of all $x \in X$ such that 
$\langle x^n : n \in B\rangle$ converges to $1$.
If $\FF$ is a filter of infinite sets, let
$\DD^X_\FF = \bigcup\{\CC^X_B : B \in \FF\}$.
The $\CC^X_B$ and $\DD^X_\FF$ are subgroups of $X$ when $X$ is abelian.
In the circle group $\TTT$, it is known that $\CC^X_B$
always has measure $0$.
We show that there is a filter $\FF$ such that
$\DD^\TTT_\FF$ has measure $0$ but is not contained in any $\CC^X_B$.
There is another filter $\GG$ such that $\DD^X_\GG = \TTT$.
We also describe the relationship between 
$\DD^\TTT_\FF$ and the 
$\DD^X_\FF$ for arbitrary compact groups $X$.
\end{abstract}

\section{Introduction} 
\label{sec-intro}
In this paper, we answer a question on topological groups asked by
Barbieri,  Dikranjan,  Milan, and Weber \cite{BDMW}, 
and we relate this question to some general facts in $C_p$
theory.  All spaces considered here are Hausdorff.
We recall the following standard definition;
see Arkhangel'skii \cite{AR} for further details on such function
spaces.

\begin{definition}
If $X,Y$ are topological spaces, then
$C(X,Y)$ is the set of continuous functions from $X$ to $Y$,
and $C_p(X,Y)$ denotes $C(X,Y)$
given the topology of \emph{pointwise convergence} (i.e., regarding
$C_p(X,Y)$ as a subset of $Y^X$ with the usual product topology).
\end{definition}
An obvious question is:
\begin{itemizz}
\item[\ding{94}]
Suppose that each $f_n\in C_p(X,Y)$
and $g$ is a limit point of
the sequence $ \langle f_n : n \in \omega\rangle$;  is there an
infinite $B \subseteq \omega$ such that the subsequence
$ \langle f_n : n \in B\rangle$ converges to $g$ pointwise
(i.e., in the topology of $C_p(X,Y)$)?
\end{itemizz}
The answer is well-known to be ``no'' in general, even when $X$ and $Y$
are both compact metric spaces. 
If $Y$ is metric, then 
for each $x\in X$, there is an infinite $B_x$ such that
$ \langle f_n(x) : n \in B_x \rangle \to g(x)$.
$B_x$ may need to vary 
with $x$, but if $X$ is compact, then one
can choose the $B_x$ from a filter; in Section \ref{sec-fs}
we show:

\begin{proposition}
\label{prop-filter}
Suppose that $X$ is compact,  $Y$ is metric, and
$g$ is a limit point of $\langle f_n : n \in \omega\rangle$ in
$C_p(X,Y)$. 
Then there is a filter $\FF \subset [\omega]^\omega$ such that 
for each $x \in X$ there is 
$B \in \FF$ such that $\langle f_n(x) : n\in B\rangle
\to g(x)$.
Furthermore, $\FF$ can be chosen to be an $F_\sigma$ subset 
of $\PP(\omega)$ (identifying $\PP(\omega)$ with $2^\omega = \{0,1\}^\omega$
in the standard way).
\end{proposition}

As usual, $\PP(\omega)$ denotes the family of all subsets of $\omega$
and $[\omega]^\omega$ denotes the family of infinite subsets of $\omega$.

One cannot replace ``metric'' by first countable here, since
the proposition fails when $X,Y$ are both the
lexicographically ordered square; see Example \ref{ex-lex}.

Observe that, in the notation of Proposition \ref{prop-filter},
$g$ must be a limit point of  $\langle f_n : n \in B\rangle$ in
$C_p(X,Y)$ for each $B \in \FF$ (see also
Lemma \ref{lemma-NN-FF}).
From this it is easily seen
that the $\FF$ is not independent of the $f_n$ and $g$.

Now, let us specialize the above discussion to
topological groups, where we consider one ``fixed'' sequences
defined uniformly across all groups:

\begin{definition}
In a group $X$, define $E_n(x) = x^n$.
\end{definition}

Note that $E_0(x) = 1$ for all $x$.
As is well-known, if $X$ is a compact topological group and
$x \in X$, then $1$ is a limit
point of the sequence $\langle x^n : n\in\omega\rangle$.
Since the same applies to every finite power $X^k$, it follows that
$E_0$ is a limit point of $\langle E_n : n \in \omega\rangle$
in $C_p(X,X)$.  In particular, if $X$ is also metric (equivalently,
second countable), then Proposition \ref{prop-filter} applies to yield
a filter.  In fact, the same filter $\FF$ can be chosen to work
for every compact metric group, since any $\FF$ which works
for the group $\TTT^\omega$ works for all compact metric groups
(see Lemma \ref{lemma-all-of}).  As usual, $\TTT$ denotes the circle group.

Returning to (\ding{94}), is there one
$B$ such that $ \langle E_n : n \in B\rangle \to E_0$?
This is answered by
Proposition \ref{prop-conv-subseq} below,
which is proved in Section \ref{sec-gps}.

\begin{definition}
If $B \in [\omega]^\omega$ and $X$ is a topological group,
let $\CC^X_B$ be the set of all $x \in X$ such that 
$\langle x^n : n \in B\rangle$ converges to $1$.
If $\FF \subset [\omega]^\omega$ is a filter, let
$\DD^X_\FF = \bigcup\{\CC^X_B : B \in \FF\}$.
\end{definition}

So, by Proposition \ref{prop-filter}, for $X$  compact metric,
one can choose $\FF$ so that $\DD^X_\FF = X$.

\begin{proposition}
\label{prop-conv-subseq}
Let $X$ be any non-trivial compact group.
If $X$ is totally disconnected, then
there is a $B \in [\omega]^\omega$ such that
$ \langle E_n : n \in B\rangle \to E_0$
(equivalently, $\CC^X_B = X$).
If $X$ is not totally disconnected, then $\CC^X_B \ne X$ 
for all $B \in [\omega]^\omega$.
\end{proposition}

Note that if $X$ is abelian, then $\CC^X_B$ and $\DD^X_\FF$
are subgroups of $X$, but this is not true in
general for non-abelian groups.
The notion of $\CC^X_B$ was considered in \cite{BDMW, CTW},
where they show that on the circle, $\CC^\TTT_B$ is a Haar nullset.
It is shown in \cite{BDMW} that assuming Martin's Axiom,
there is a Haar null subgroup $N$ of $\TTT$ which is not contained
in any $\CC^\TTT_B$, and they raise the question of whether
Martin's Axiom was needed.  Here, we produce such an $N$
without the use of any set-theoretic axioms.  More generally,

\begin{theorem}
\label{thm-nicefilter}
Let $\FF \subset [\omega]^\omega$ be the filter generated by all
sets of the form $\{k! + 1 : k \in D\}$, where $D\subseteq \omega$
has asymptotic density $1$.  Then:
\begin{itemizz}
\item[1.] $\FF$ is a Borel subset of $\PP(\omega) \cong 2^\omega$.
\item[2.] Whenever $X$ is an infinite compact group:
\begin{itemizz}
\item[a.] $\DD^X_\FF$ is a Haar nullset.
\item[b.] If $X$ is not totally disconnected, then
$\DD^X_\FF$ is not a subset of $\CC^X_B$ for any infinite $B$.
\end{itemizz}
\end{itemizz}
\end{theorem}
In particular, when $X$ is abelian and not totally disconnected,
$\DD^X_\FF$ is a Haar 
null subgroup not contained in any $\CC^X_B$.
Of course, (1) of this theorem is obvious; (2) is proved in 
Section \ref{sec-gps}.
Observe, by Proposition \ref{prop-conv-subseq},
that we cannot delete the ``not totally disconnected'' in (2b).

Section \ref{sec-gps} contains some further information
about properties of the $\CC^X_B$.
In many cases, these properties for an arbitrary
$X$ can be inferred directly from the special cases $X = \TTT$
or $X = \TTT^\omega$.

We conclude this Introduction with some easy remarks
which will simplify the notation in the next two sections.

\begin{definition}
\label{def-tilde}
If $\FF \subset [\omega]^\omega$ is a filter, then
$\widetilde \FF$ is the filter generated by $\FF$ and all cofinite
sets.
If $B \in [\omega]^\omega$ then $B \up$ is the filter
$\{W \subseteq \omega: B \subseteq W \}$.
\end{definition}

So, $\widetilde{B \up\,} = \{W\subseteq \omega : B \subseteq^* W \}$.
Since the convergence of a sequence
$ \langle f_n(x) : n \in B\rangle$ does not change if one modifies
$B$ on a finite set, we can always deal with $\widetilde \FF$
rather than $\FF$.  In particular,

\begin{lemma}
\label{lemma-same-tilde}
$\DD^X_\FF = \DD^X_{\widetilde \FF}$.
\end{lemma}

Any general lemma about all $\DD^X_\FF$ will also apply to
the $\CC^X_B$ by:

\begin{lemma}
\label{lemma-set-filter}
$\CC^X_B = \DD^X_{B \up}$.
\end{lemma}

\section{Filters}
\label{sec-fs}

\begin{proofof}{Proposition \ref{prop-filter}}
Let $d$ be a metric on $Y$.  Choose (inductively) disjoint finite
$S_k \subset \omega$ such that $S_0 = \emptyset$, and for $k\ge 1$,
$$
\forall x_0, \ldots x_{k-1}\in X\;  \exists n \in S_k \;
\forall i < k\;  [ d(f_n(x_i), g(x_i)) < 1/k]\ \ .  \eqno{(\mbox{\ding{76}})}
$$
To see that such a finite $S_k$ exists, let
$A_k = \omega \setminus \bigcup_{j<k} S_j$, and let
$$
U_{n,k} = \{( x_0, \ldots x_{k-1})\in X^k :
\forall i < k\;  [ d(f_n(x_i), g(x_i)) < 1/k]\}  \ \ .
$$
Then $\bigcup\{U_{n,k} : n \in A_k\} = X^k$ (since 
$g$ is a limit point of $\langle f_n : n \in A_k \rangle$),
so by compactness of $X^k$, there is a finite $S_k \subset A_k$
such that $\bigcup\{U_{n,k} : n \in S_k\} = X^k$.

For $x\in X$, let 
$B_x = \bigcup_{k \in \omega}\{n \in S_{k} : d(f_n(x),g(x)) < 1/k \}$.
Applying (\mbox{\ding{76}}), each 
$B_{x_0} \cap \cdots \cap  B_{x_{\ell-1}}$ meets $S_k$ for all $k \ge \ell$,
and is hence an infinite set.
Let $\FF$ be the filter generated by $\{B_x : x \in X\}$.
This $\FF$ satisfies the theorem because 
$\langle f_n(x) : n\in B_x \rangle \to g(x)$.

To see that $\FF$ is an $F_\sigma$ in $\PP(\omega) \cong 2^\omega$, let
$$
\GG_\ell = \{D \in \PP(\omega) : 
\exists x_0, \ldots x_{\ell-1}\in X\;  
[B_{x_0} \cap \cdots \cap  B_{x_{\ell-1}} \subseteq D] \} \ \ .
$$
Then  $\FF = \bigcup_\ell \GG_\ell$.  To prove that $\GG_\ell$ is closed,
let
$$
\begin{array}{l}
\HH_\ell = \left\{(x_0, \ldots, x_{\ell-1}, D) \in X^\ell \times \PP(\omega)
\ : \ \forall k \in \omega \forall n \in S_k \right.  \\
\left.
\left[ d(f_n(x_0), g(x_0)) < 1/k  \;\wedge\; \ldots \;\wedge\;
d(f_n(x_{\ell-1}), g(x_{\ell-1})) < 1/k  
\  \Longrightarrow \  n \in D \right] \right\} \ .
\end{array}
$$
Then $\HH_\ell$ is closed in $X^\ell \times \PP(\omega)$.  Hence,
$\GG_\ell \subseteq \PP(\omega)$, which is the projection of $\HH_\ell$,
is closed as well.
\end{proofof}

The idea for obtaining the $S_k$ above is taken from the proof that
$C_p(X,Y)$ has countable tightness
(Kelley and Namioka \cite{KeNa}, Lemma 8.19; see also \cite{AR}, Ch.~II\S1).
We do not have, in general, a simpler description of a filter 
$\FF$ satisfying
Proposition \ref{prop-filter}, although $\FF$ is
related to the neighborhood filter in $C_p(X,Y)$:

\begin{definition}
\label{def-neigh}
If $g$ is a limit point of $\langle f_n : n \in \omega\rangle$ in
$C_p(X,Y)$, then the \emph{induced neighborhood filter}
is the filter generated by all subsets of $\omega$ of the form
$\{n\in\omega : f_n \in U\}$, where $U$ is a neighborhood
of $g$ in $C_p(X,Y)$.
\end{definition}

\begin{lemma}
\label{lemma-NN-FF}
If $\FF$ is as in Proposition \ref{prop-filter} and $\NN$
is the induced neighborhood filter, then $\NN \subseteq \widetilde \FF$
(see Definition  \ref{def-tilde}).
\end{lemma}
\begin{proof}
$\NN$ is generated by sets of the form
$A = \{n\in\omega: f_n(x) \in V\}$, where
$x \in X$ and
$V$ is a neighborhood of $g(x)$ in 
$Y$.  Fix $B\in\FF$ with 
$\langle f_n(x) : n \in B\rangle$ converging to $g(x)$.
Then $B\subseteq^* A$, so $A\in \widetilde\FF$.
\end{proof}

In some cases,
one can simply take $\FF = \NN$ in Proposition \ref{prop-filter};
for example, this will always work if $Y$ is finite.  
However, this does not work in the case $X = Y = \TTT$, with
$f_n = E_n$ and $g = E_0$; see Proposition \ref{proposition-NN-rat}.

\section{Groups}
\label{sec-gps}
We consider now in more detail the
$\CC^X_B$ and $\DD^X_\FF$ for compact groups $X$.
Basic facts about such groups can be found in Hofmann and Morris \cite{HM}.

\begin{lemma}
\label{lemma-0-dim}
If $X$ is a totally disconnected compact group,
and $B = \{k! : k \in \omega\}$, then $\CC^X_B = X$.
\end{lemma}
\begin{proof}
Such an $X$ is an inverse limit of finite groups,
or, equivalently,
a closed subgroup of a product of finite groups
(see \cite{HM}, Theorem 1.34).
\end{proof}

For the other groups, we shall show that
$\CC^X_B$ is never all of $X$, while
$\DD^X_\FF$ may or may not be all of $X$, depending on $\FF$ and $X$.

First, consider the ``all of $X$'' case:

\begin{lemma}
\label{lemma-all-of}
If $\DD^{\TTT^\omega}_\FF = {\TTT^\omega}$ then $\DD^X_\FF = X$ for
all compact metric groups $X$.
\end{lemma}
\begin{proof}
Since each $x \in X$ generates
a compact abelian subgroup,
we may assume that $X$ is abelian, in which case $X$ is continuously
isomorphic to a subgroup of $\TTT^\omega$.
\end{proof}

For the ``not all of $X$'' case, there is a similar reduction to $\TTT$.
We shall show:

\begin{theorem}
\label{thm-reduce}
Let $\FF \subset [\omega]^\omega$ be a filter which is analytic
as a subset of $\PP(\omega) \cong 2^\omega$, and 
assume that $\DD^\TTT_\FF \ne \TTT$.
Let $X$ be any infinite compact group.  Then
$\DD^X_\FF$ is a Haar nullset if at least one of the following holds:
\begin{itemizz}
\item[1.] $X$ is abelian and not totally disconnected.
\item[2.] $X$ is connected.
\item[3.] For all $B \in \FF$ and all $m \ge 2$:
$\{n \in B : m \nmid n \}$ is infinite.
\end{itemizz}
Furthermore, $\DD^X_\FF \ne X$ whenever $X$ is not totally disconnected.
\end{theorem}

This in particular applies when $\FF = B\up$
(see Definition \ref{def-tilde}
and Lemmas \ref{lemma-same-tilde} and \ref{lemma-set-filter}).
Observe that

\begin{lemma}
\label{lemma-T-null}
$\CC^\TTT_B$ is a Haar nullset for all infinite $B \subseteq \omega$.
\end{lemma}

See \cite{CTW} (Lemma 3.10) and \cite{BDMW} \S4 for proofs;
yet another proof is given below; see Remark \ref{remark-null-CC}.
It follows that $\CC^X_B$ is a Haar nullset if one of (1,2,3)
from Theorem \ref{thm-reduce} hold.
Note that these conditions cannot be dropped.  For example,
let $X = O(2)$.  Then the component of $1$ in $X$ is
$SO(2) \cong \TTT$, which has index $2$ in $X$.  
The two cosets of $SO(2)$ are $SO(2)$ and $R$ (the reflections),
and the elements of $R$  all have order $2$.  
Now, let $B = \{2k : k \in \omega\}$.  Then $\CC^X_B$ contains
all of $R$ plus two elements of $SO(2)$, so $\CC^X_B$
has Haar measure $1/2$.

Observe that if $\FF$ is analytic, then $\DD^\TTT_\FF$ is also analytic,
and hence Haar measurable.
Since a measurable subgroup is null iff it has infinite index,
the following lemma, called the Steinhaus-Weil Theorem in \cite{CTW},
is relevant:

\begin{lemma}
\label{lemma-finite-index}
If $X$ is a compact group and $H$ is a measurable subgroup of
finite index, then $H$ is clopen in $X$.
\end{lemma}
\begin{proof}
$H$ has positive measure, so $H\iv H = H$ has non-empty interior
(see \cite{HR}, Cor.~20.17),
so $H$ is open (since it is a group), and hence clopen
(since it has finite index).
\end{proof}

In particular, if $X$ is connected, then either $H$ is null or $H = X$.  Thus,

\begin{corollary}
\label{cor-null-T}
If $\FF$ is an analytic filter and 
$\DD^\TTT_\FF \ne \TTT$, then
$\DD^\TTT_\FF$ is a Haar nullset.
\end{corollary}

Now, to prove Theorem \ref{thm-reduce},
we prove a few lemmas which reduce the situation for a general $X$ 
to the case $X = \TTT$.

\begin{lemma}
\label{lemma-map-onto}
If $X,Y$ are compact groups and $\varphi$ is a continuous
homomorphism from $X$ onto $Y$, then
$\DD^X_\FF \subseteq \varphi\iv (\DD^Y_\FF)$.
\end{lemma}

Since $\varphi\iv$ also preserves Haar measure, one may prove
that $\DD^X_\FF$ is null by proving 
that $\DD^Y_\FF$ is null.  Using this remark,
the reduction for abelian groups is easy:

\begin{lemma}
\label{lemma-reduce-to-T}
If  $\DD^\TTT_\FF$ is null, then
$\DD^X_\FF$ is null for all infinite compact 
abelian groups $X$ which are not totally disconnected.
\end{lemma}
\begin{proof}
$X$ is continuously isomorphic to some
compact subgroup of $\TTT^\theta$ for some cardinal $\theta$,
so we may assume
that $X \subseteq \TTT^\theta$.
If $\pi_\alpha$ is the projection onto the $\alpha^{\mathrm{th}}$
coordinate, then $\pi_\alpha(X)$ is a compact subgroup of $\TTT$.
Some of the $\pi_\alpha(X)$ may be finite, but they 
cannot all be finite unless $X$ is totally disconnected.
However, if $\pi_\alpha(X)$ is infinite, then it must be all of $\TTT$,
so the lemma follows by using Lemma \ref{lemma-map-onto}
with $\varphi = \pi_\alpha$.
\end{proof}

To handle non-abelian $X$, we shall replace $\TTT$ in the above
proof by a compact Lie group.
For this paper, we can take as a definition
that $X$ is a \textit{compact Lie group} iff $X$ is continuously 
isomorphic to a compact subgroup of the unitary group
$U(n)$ for some finite $n$.
Many other equivalents are known; see \cite{HM, KN}.
Observe that all finite groups are compact Lie groups by this definition.

\begin{lemma}
\label{lemma-reduce-to-lie}
If $X$ is a compact group and is not totally disconnected,
then there is a continuous homomorphism $\pi$ from $X$ onto
some infinite compact Lie group $Y$.
\end{lemma}
\begin{proof}
By standard representation theory (see \cite{HR, HM}),
we may assume that
$X \subseteq \prod_{\alpha  < \theta} U(n_\alpha)$.
Then each $\pi_\alpha(X)$ is a Lie group, and
at least one of the $\pi_\alpha(X)$ is infinite.
\end{proof}

Then, as in Lemma \ref{lemma-reduce-to-T}, if
$\DD^Y_\FF$ is null then $\DD^X_\FF$ is null.
However, in proving that $\DD^Y_\FF$  is null from the assumption
that $\DD^\TTT_\FF$ is null, we cannot use a similar argument, 
since $Y$ need need not have a homomorphism onto $\TTT$.
Rather, we use the fact that $\TTT$ is contained in $Y$.
We shall apply Lemma \ref{lemma-reduce-to-T} to the maximal
tori in $Y$ (i.e., maximal connected abelian subgroups; see \cite{HM,KN})
to get:

\begin{lemma}
\label{lemma-conn-lie}
If  $\DD^\TTT_\FF$ is null, then
$\DD^Y_\FF$ is null for all non-trivial connected compact Lie groups $Y$.
\end{lemma}
\begin{proof}
Let $f$ be the characteristic function of $\DD^Y_\FF$;
so we wish to show that $\int_Y f \, d \lambda_Y = 0$,
where $\lambda_Y$ is normalized Haar measure on $Y$.
By Lemma \ref{lemma-reduce-to-T}, 
we know that $\int_H f \, d \lambda_H = 0$ whenever $H$
is a maximal torus in $Y$.
Now, observe that $f$ is a class function; that is
$f(x\iv y x) = f(y)$.
It follows that we may integrate $f$ by the Weyl Integration Formula
(see \cite{KN}, eqn.~(8.62)):
$$
\int_Y f(x) \, d \lambda_Y(x) = {1 \over |W|}
\int_H f(t)|D(t)|^2 \, d \lambda_H(t) \ \ ,
$$
where $W$ is a finite group and $D(t)$ is a finite function of $t$.
Since $f(t) = 0$ for $\lambda_H$ -- almost every $t$, we get
$\int_Y f \, d \lambda_Y = 0$.
\end{proof}

Note that $Y$ must be connected for this 
Weyl Integration Formula
to be true, since all the maximal tori are contained 
in the identity component of $Y$; furthermore,
as pointed out above, the lemma fails for $Y = O(2)$.

\begin{proofof}{Theorem \ref{thm-reduce}}
$\DD^\TTT_\FF$ is null by Corollary \ref{cor-null-T}.
Now, assume that $X$ is not totally disconnected.
Then $\DD^X_\FF$ is null if $X$ is abelian by Lemma \ref{lemma-reduce-to-T},
so we must handle the non-abelian case.
If $X$ is connected, then we can map $X$ onto a non-trivial connected
compact Lie group $Y$ by Lemma \ref{lemma-reduce-to-lie}, so that 
$\DD^X_\FF$ is null by Lemmas 
\ref{lemma-conn-lie} and \ref{lemma-map-onto}.
If $X$ is not connected, then $Y$ may fail to be connected.
If $Y_0$ is the identity component of $Y$, then
$\DD^Y_\FF\cap Y_0 = \DD^{Y_0}_\FF$ will still
be null by Lemma \ref{lemma-conn-lie},
which proves that  $\DD^X_\FF \ne X$. 

Now, we must prove that $\DD^X_\FF$ is null in Case (3) of
Theorem \ref{thm-reduce}.
First, assume that $X$ is not totally disconnected.
Using the same $Y$, it is
again sufficient to show that $\DD^Y_\FF$ is null.
But we already know that $\DD^Y_\FF\cap Y_0$ is null, and 
in Case (3), $\DD^Y_\FF \subseteq Y_0$.  To see this,
fix $y \in Y \backslash Y_0$, and let $m$ be the order
of $[y]$ in the quotient $Y/Y_0$.
If $y$ were in $\DD^Y_\FF$, we could fix $B \in \FF$ such that 
$\langle y^n : n \in B\rangle \to 1$.
But $C = \{n \in B : m \nmid n \}$ is infinite,
$Y_0$ is a neighborhood of $1$, and 
$y^n\notin Y_0$ for all $n \in C$, a contradiction.

Finally if $X$ is totally disconnected in Case (3),
then $\DD^X_\FF = \{1\}$, as in the proof of Lemma \ref{lemma-0-dim}.
\end{proofof}

\begin{proofof}{Proposition \ref{prop-conv-subseq}}
By Lemmas \ref{lemma-0-dim} and \ref{lemma-T-null}
and Theorem \ref{thm-reduce}.
\end{proofof}

We now give some examples of $\FF$ for which $\DD^\TTT_\FF$ is null,
using some well-known facts about Hadamard sets and
the Bohr topology on $\ZZZ$;
see \cite{KR} for definitions and references
to the earlier literature.

\begin{definition}
\label{def-bohr-N}
$\NNb\subseteq \PP(\omega)$ denotes the neighborhood
filter at $0$ in the topology $\omega$ inherits as
a subset of $\ZZZ^\#$ (that is, the group $\ZZZ$ with
its Bohr topology).
\end{definition}

Equivalently, $\NNb$ is the induced neighborhood filter
in the sense of Definition \ref{def-neigh},
taking $f_n = E_n$ and $g = E_0$.
Applying Lemma \ref{lemma-NN-FF},

\begin{lemma}
\label{lemma-NF}
If $\DD^\TTT_\FF = \TTT$ then $\NNb \subseteq \widetilde\FF$.
\end{lemma}

Note that $\DD^\TTT_\NNb$ is countable
(see Proposition \ref{proposition-NN-rat}),
so to get $\DD^\TTT_\FF = \TTT$
(as in Proposition \ref{prop-filter}), $\FF$ must properly extend $\NNb$.

Applying Lemma \ref{lemma-NF} and Corollary \ref{cor-null-T},

\begin{lemma}
\label{lemma-null-DD}
Suppose that 
$\NNb \not\subseteq \widetilde\FF$, and $\FF$ is analytic as a subset
of $\PP(\omega) \cong 2^\omega$.
Then $\DD^\TTT_\FF$ is a Haar nullset.
\end{lemma}

\begin{remark}
\label{remark-null-CC}
{\rm
This yields another proof of Lemma \ref{lemma-T-null}:
apply Lemma \ref{lemma-null-DD} to
$\FF = B\up$, and note that 
$\NNb \not\subseteq \widetilde\FF$ because $\ZZZ^\#$ 
has no convergent $\omega$-sequences.
}
\end{remark}

\begin{proofof}{Theorem \ref{thm-nicefilter}(2a)}
$H =  \{k! +1 : k\in\omega\}$ is a Hadamard set, and
hence closed and discrete in $\ZZZ^\#$.
Thus, $\omega \backslash H \in \NNb$, so that $\NNb \not\subseteq \FF$.
Then, (2a) follows by Lemma \ref{lemma-null-DD} and
Theorem \ref{thm-reduce} (Case (3)).
\end{proofof}

Now, let us turn to a proof of Theorem \ref{thm-nicefilter}(2b).
Here, most of the work will be done on the \textit{solenoid}
$\widehat \QQQ$ (see \cite{HR}):

\begin{definition}
\label{def-qhat}
$\widehat \QQQ$ denotes the dual of the discrete group of rationals.
We shall realize $\widehat \QQQ$ concretely as
$\{\vec z \in \TTT^\omega :
\forall \alpha < \omega [ (z_{\alpha+1})^{\alpha+1} = z_\alpha ]\}$.
\end{definition}

This differs slightly from the notation in \cite{HR}.
We can identify $\vec z$ with the character of $\QQQ$ which
takes  $1/\alpha!$ to $z_\alpha$.

\begin{lemma}
\label{lemma-get-solenoid}
If  $X$ is an infinite compact group and is not totally disconnected,
then $X$ has a non-trivial closed subgroup which is
a continuous homomorphic image of $\widehat \QQQ$.
\end{lemma}
\begin{proof}
First, WLOG, $X$ is abelian.  To see this, let
$\pi : X \to Y$ be a continuous homomorphism onto an infinite compact
Lie group $Y$ (see Lemma \ref{lemma-reduce-to-lie}).
Choose $x \in X$ such that $y := \pi(x)$ has infinite order.  Then $\pi$ maps
$\overline {\langle x \rangle }$ onto
$\overline {\langle y \rangle }$, which is an infinite compact
Lie group, and hence not totally disconnected.  It follows
that $\overline {\langle x \rangle }$ is not totally disconnected
(see \cite{HM}, Exercise E1.13),
so we may replace $X$ by $\overline {\langle x \rangle }$.

Second, WLOG, $X$ is connected, since we may replace $X$ by
the component of $1$.

Now, let $G = \widehat X$, which is a discrete torsion-free
abelian group.  Let $V \supseteq G$ be the divisible hull of $G$;
then $V$ is torsion-free and divisible, so we may regard $V$
as a vector space over $\QQQ$.   Choose a basis
$\{v_\alpha : \alpha < \kappa\}$ for $V$ with $v_0 \in G$,
let $W$ be the vector subspace generated by $v_0$,
and let $\psi$ be the canonical homomorphism of $V$ onto $W$.
Let $H = \psi(G) \subseteq W$.  Then $H$ is torsion-free
and non-trivial (since $v_0 \in G$), and $\psi\res G$ maps
$G$ onto $H$, so $Y := \widehat{H}$ is a non-trivial closed
subgroup of $X = \widehat{G}$.
But also, $H$ is isomorphic to a subgroup of $\QQQ \cong W$,
so $Y$ is a quotient of $\widehat \QQQ$.
\end{proof}

\begin{definition}
\label{def-thin}
$A \subseteq \omega$ is \emph{thin} iff $A$ is infinite and
of the form $\{a_k : k\in \omega\}$, where
$0 < a_0 < a_1 < \cdots$ and $\lim_k a_k/a_{k+1} = 0$.
\end{definition}

\begin{lemma}
\label{lemma-limits-in-solenoids}
If
$A \subseteq \omega$ is thin,
$A$ is partitioned into disjoint infinite subsets $B,C$,
and $\vec v, \vec w \in \widehat \QQQ$, then there is a
$\vec z \in \widehat \QQQ$ such that
$\langle (\vec z)^n : n \in B\rangle$ converges to $\vec v$ and
$\langle (\vec z)^n : n \in C\rangle$ converges to $\vec w$.
\end{lemma}
\begin{proof}
List $A$ in increasing order as $\{a_j : j\in \omega\}$.
Let $\varepsilon_j = \sup_{\ell \ge j} (a_\ell / a_{\ell+1})$.
Then $\varepsilon_j \searrow 0$, each $a_j / a_{j+1} \le \varepsilon_j$,
and $\varepsilon_j \le 1$.

Let $\gamma_0 = 0$, and let $\gamma_{j+1}$ be
the largest integer $\gamma$ such that $\gamma ! \le 1/\sqrt{\varepsilon_j}$.
Note that $\gamma_j \nearrow \infty$ and
$$
(a_j/a_{j+1}) (\gamma_{j+1} ! )\le \varepsilon_j (\gamma_{j+1}!) \le
\sqrt{\varepsilon_j} \ \ .
$$
Hence,
$(a_j/a_{j+2}) (\gamma_{j+2} ! )\le
(a_j/a_{j+1}) (\gamma_{j+1} ! ) \cdot (a_{j+1}/a_{j+2}) (\gamma_{j+2} ! )\le
\sqrt{\varepsilon_j \varepsilon_{j+1}}$.
Continuing in this way,
$$
j < k \ \Longrightarrow\   (a_j/a_{k}) (\gamma_{k}! )\le
\sqrt{\varepsilon_j \varepsilon_{j+1} \cdots \varepsilon_{k-1}}
\le (\sqrt{\varepsilon_j})^{k-j} \ \ .
\eqno{(\mstar)}
$$

Our $\vec z$ will be $\lim_j \vec z_j$, where each $\vec z_j\in \widehat \QQQ$;
the $\vec z_j$ will be chosen by induction on $j$.
Use $z_{j,\alpha}$ for the components of
$\vec z_j \in \widehat\QQQ \subset \TTT^\omega$.
Let $d$ be the metric on $\TTT$ obtained by
identifying $\TTT$ with $\RRR/\ZZZ$, so that the circumference of $\TTT$
is $1$, and $d(e^{i\theta}, 1) = |\theta|/2\pi$ when $|\theta| \le \pi$.
We shall get, for each $k \in \omega$:
\begin{itemizz}
\item[1.] For $\alpha \le \gamma_k$:
$(z_{k,\alpha})^{a_k}$ equals $v_\alpha$ when $a_k \in B$ and
$w_\alpha$ when $a_k \in C$.
\item[2.]For $\alpha \le \gamma_{k}$, when $k>0$:
$d(z_{k, \alpha}, z_{k-1, \alpha}) \le (\gamma_{k}!)/a_{k}$.
\end{itemizz}
To ensure (1) for a given $k$, it is sufficient to choose
$z_{k,\gamma_k}$ so that (1) holds for $\alpha = \gamma_k$;
then (1) will hold for $\alpha < \gamma_k$ by our definition
(\ref{def-qhat}) of $\widehat \QQQ$.
For $\alpha > \gamma_k$, we just choose the $z_{k,\alpha}$ so that
the point $\vec z_k$ lies in $\widehat \QQQ$.
To get (2) along with (1):  We are given $\vec z_{k-1}$ and we must
define $\vec z_{k}$ by choosing $z_{k,\gamma_{k}}$
so that $(z_{k,\gamma_{k}})^{a_{k}} $ is $v_{\gamma_{k}}$
when $a_{k} \in B$ and
$w_{\gamma_{k}}$ when $a_{k} \in C$.
There are $a_{k}$ possible choices for $z_{k,\gamma_{k}}$,
spaced evenly around the circle, at distance $1 / a_{k}$ apart,
so we can make the choice so that
$d(z_{k, \gamma_{k}}, z_{{k-1}, \gamma_{k}}) \le 1/ a_{k}$.
Then (2) follows by our definition of $\widehat \QQQ$.

Now, if we fix $j$ and $\alpha \le \gamma_j$, then
whenever $j < k$, we may apply (2) and (\mstar) to get
$a_j d(z_{k, \alpha}, z_{k-1, \alpha}) \le
(a_j /a_k) (\gamma_{k}!) \le (\sqrt{\varepsilon_j})^{k-j}$.
Thus, whenever $\ell > j$:
$$
a_j d(z_{\ell, \alpha}, z_{j, \alpha}) \le
\sum_{k=j+1}^\ell (\sqrt{\varepsilon_j})^{k-j} \le 
\sqrt{\varepsilon_j}/(1 - \sqrt{\varepsilon_j})  \ \ .
\eqno{(\mdag)}
$$
In particular, the sequence $\langle z_{j, \alpha} : \alpha \in \omega \rangle$
is Cauchy, so we can define $\vec z = \lim_j \vec z_j$.
Then, (\mdag) yields
$d(z_{\alpha}^{a_j}, z_{j, \alpha}^{a_j}) \le
\sqrt{\varepsilon_j}/(1 - \sqrt{\varepsilon_j})$
whenever $\alpha \le \gamma_j$.
Applying this for the $a_j \in B$, when
$(z_{j,\alpha})^{a_j} = v_\alpha$, we get
$\langle (\vec z)^n : n \in B\rangle \to \vec v$, and
applying it for the $a_j \in C$ yields
$\langle (\vec z)^n : n \in C\rangle \to \vec w$.
\end{proof}

\begin{lemma}
\label{lemma-thin-limits}
Let $X$ be any infinite compact group which is not totally disconnected.
Let $A \subseteq \omega$ be thin, with
$A$ partitioned into disjoint infinite subsets $B,C$.
Then there are $x,y \in X$ with $y \ne 1$ such that
$\langle x^n : n \in B\rangle$ converges to $y$ and
$\langle x^n : n \in C\rangle$ converges to $1$.
\end{lemma}
\begin{proof}
By Lemma \ref{lemma-get-solenoid}, we may assume that
there is a continuous homomorphism $\varphi$ from $\widehat \QQQ$
onto $X$.  Then $y$ can be any element of $X \backslash \{1\}$.
Now, in $\widehat \QQQ$,  let $\vec w = 1 $, and choose
$\vec v$ with $\varphi(\vec v) = y$, and apply
Lemma \ref{lemma-limits-in-solenoids}.
\end{proof}

Consider this in particular with $X = \TTT$, where the
argument of Lemma \ref{lemma-limits-in-solenoids} can be
done (considerably simplified) in $\TTT$ directly, and
resembles the proof that Hadamard sets 
are $I_0$ sets.  However, Lemma \ref{lemma-thin-limits} fails
if we only assume that $A$ is a Hadamard set.
For example, let $A = \{2^j : j\in \omega\}$,
$B = \{2^{2k+1} : k \in \omega\}$, and
$C = \{2^{2k} : k \in \omega\}$, and consider $x\in\TTT$.
If $\langle x^n : n \in C\rangle \to 1$, then
$x^{2^j} = 1$ for some $j$, but then also
$\langle x^n : n \in B\rangle \to 1$.

\begin{proofof}{Theorem \ref{thm-nicefilter}(2b)}
Fix an infinite $B$; we need to show that
$\DD^X_\FF \not\subseteq \CC^X_B$.
Since $\CC^X_B$ gets bigger as $B$ gets smaller, we may assume
that $B$ is small enough so that for some $C \in \FF$:
$C \cap B = \emptyset$ and $A := C \cup B$ is thin.
Then, applying Lemma \ref{lemma-thin-limits}, we may fix $x\in X$
so that $\langle x^n : n \in C\rangle$ converges to $1$
and $\langle x^n : n \in B\rangle$ converges to $y \ne 1$.
Then $x \in \DD^X_\FF$ and $x \notin \CC^X_B$.
\end{proofof}

\section{Examples}
\label{sec-ex}
First, we point out that Proposition \ref{prop-filter} can 
fail if $X$ is not assumed to be compact.

\begin{example}
\label{ex-non-p}
Let $X$ be the cardinal $\cccc = 2^{\aleph_0}$ with the discrete topology,
and let $Y = \omega + 1$ with the order topology.
Then in $C_p(X,Y)$, there is a
sequence $\langle f_n : n \in \omega\rangle$ with limit point $g$
such that there is no filter
$\FF \subset [\omega]^\omega$ satisfying:
$$
\forall x \in X \exists B \in \FF [ \langle f_n(x) : n\in B\rangle \to g(x) ]
\eqno{(\mbox{\ding{40}})}
$$
\end{example}
\begin{proof}
Let $\UU$ be a non-principal ultrafilter on $\omega$ which 
is not a P-point.  List $\UU$ as $\{B_\alpha : 0 < \alpha < \cccc\}$.
Also, partition $\omega$ into infinite sets $A_k$ (for $k \in \omega$)
so that each $A_k \notin \UU$ and each $B \in \UU$ meets
some $A_k$ in an infinite set.

Let $g(\alpha) = \omega$ for all $\alpha$.
For $\alpha > 0$, let $f_n(\alpha)$ be $n$ for $n \in B_\alpha$ and
$0$ for $n \notin B_\alpha$.  Let $f_n(0)$ be the $k$ such 
that $n \in A_k$.

Now, suppose that $\FF$ satisfied (\ding{40}).
Fix $C \in \FF$  with $\langle f_n(0) : n\in C\rangle \to \omega$.
Then each $C \cap A_k$ is finite, so $C \notin \UU$, so fix 
$\alpha > 0$ with $B_\alpha \cap C = \emptyset$.
Now fix $D \in \FF$ with $\langle f_n(\alpha) : n\in D\rangle \to \omega$.
Then $D \subseteq^* B_\alpha$, so $D \cap C$ is finite.
But $D \cap C \in \FF$, contradicting $\FF \subset [\omega]^\omega$.
\end{proof}

This is the simplest possible $Y$ for such an example, since
Proposition \ref{prop-filter} does hold for arbitrary $X$
whenever $Y$ is finite,
taking $\FF$ to be the induced neighborhood filter $\NN$
(as in Definition \ref{def-neigh}).
Also, under Martin's Axiom, Proposition \ref{prop-filter} holds
for all $X$ of size less than $\cccc$ whenever $Y$ is first countable.
When $|X| = \cccc$,
taking $Y$ to be first countable is not enough,
even when $X$ compact, since Proposition \ref{prop-filter} fails
when $X$ and $Y$ are the lexicographically ordered square.
This is the space $[0,1]\times[0,1]$, ordered lexicographically,
and given the usual order topology; note that it is compact
and first countable.

\begin{example}
\label{ex-lex}
Let $X$ be the  lexicographically ordered square.
Then there is a
sequence $\langle f_n : n \in \omega\rangle$ in $C_p(X,X)$
with limit point $g$
such that there is no filter
$\FF \subset [\omega]^\omega$ satisfying (\ding{40}) above.
\end{example}
\begin{proof}
Let
$\langle \widetilde{f_n}: n \in \omega\rangle$  and $\widetilde g$ 
be the functions in $C_p([0,1]_d, \,\omega+1)$ obtained
in Example \ref{ex-non-p}, identifying the ordinal $\cccc$ with
the discrete $[0,1]$.
We shall encode this example into $C_p(X,X)$.

To embed $\omega + 1$ into $[0,1]$,  
let $u_\omega = 1/2$ and fix $u_n$ with $0 < u_n < 1/2$
and $u_n \nearrow u_\omega$.
Next, in $X$, let $I_r = \{r\}\times [0,1]$; then $I_r$ is
a homeomorphic copy of $[0,1]$ for each $r \in [0,1]$.
We shall have $f_n : X \to X$ map each $I_r$ onto $I_r$,
and encode the $\widetilde{f_n}$ as follows:
Let
$$
f_n(r,1/2) = (r,u_{\widetilde{f_n}(r)}) 
\ , \ f_n(r,0) = (r,0)
\  , \ f_n(r,1) = (r,1) \ \ .
$$
Then, fill in the rest of the values $f_n(r,y)$ for $y \in [0,1]$
by linear interpolation, so that the graph of $f_n$ intersected with
$I_r \times I_r$ is the union of two line segments.

Likewise, we get $g$ from $\widetilde g$,
but then $g$ is the identity function.
We can now apply the same argument as in Example \ref{ex-non-p}.
\end{proof}

We remark on the relevance of P-points here.
For arbitrary $X$ and first countable $Y$, 
let $\NN$ be the induced neighborhood filter
(as in Definition \ref{def-neigh}).
If there is a P-point $\UU$ such that
$\NN$ is contained in $\UU$,
then taking $\FF = \UU$ will work in Proposition \ref{prop-filter}.
In Examples \ref{ex-non-p} and \ref{ex-lex}
we have $\NN = \UU$, but that $\UU$ was 
chosen to be a \textit{non}-P-point. 
For $X$ compact and $Y$ metric,
if we assume Martin's Axiom, then we can produce a P-point ultrafilter
$\FF \supseteq \NN$ which gives the convergence described in
Proposition \ref{prop-filter}, but, of course,
this $\FF$ will not be a Borel set.
The proof of Proposition \ref{prop-filter} amounts to adding to 
$\NN$ some diagonal intersections of sets in $\NN$, which
is part of the P-point construction, although we do not actually
build an ultrafilter.
Note that $\NN$ itself will not in general satisfy
Proposition \ref{prop-filter}; for example,
to get a filter $\FF$ with $\DD^\TTT_\FF = \TTT$, we cannot just take 
$\FF = \NNb$ (see Definition  \ref{def-bohr-N}):

\begin{proposition}
\label{proposition-NN-rat}
$\DD^\TTT_\NNb$ is the set of all $w \in \TTT$
of finite order.
\end{proposition}
\begin{proof}
If $w^k = 1$, then $w^n = 1$ for all $n \in \{jk : j \in \omega \} \in \NNb$.

Conversely, suppose $w$ has infinite order.
Consider $\ZZZ \subseteq \bohr \ZZZ = \Hom(\TTT_d, \ZZZ)$ in the standard way,
and let $\lambda$ be Haar measure on $\bohr \ZZZ$.
Observe that if $B \in \NNb$, then $\lambda (\overline B) > 0$,
whereas if $ \langle w^n : n \in B\rangle \to 1$,
then $\lambda (\overline B) = 0$.
\end{proof}

The following shows that Lemma \ref{lemma-all-of}
can fail if $X$ is non-metrizable:

\begin{example}
\label{example-T-c}
If $\DD^{\TTT^\omega}_\FF = \TTT^\omega$, then
$\DD^{\TTT^\cccc}_\FF$ has inner Haar measure $0$
and outer Haar measure $1$.
\end{example}
\begin{proof}
Let $P$ be the set of all $f \in {\TTT^\cccc}$ such 
that $f$ maps $\cccc$ onto $\TTT$, and let
$Q$ be the set of elements of $f \in {\TTT^\cccc}$
such that $\{\alpha : f(\alpha) \ne 1\}$ is countable.
Then $Q$ has outer measure $1$ and $Q \subseteq \DD^{\TTT^\cccc}_\FF$
(since $\DD^{\TTT^\omega}  = \TTT^\omega)$, so 
$\DD^{\TTT^\cccc}_\FF$ has outer measure $1$.
Also, $P$ has outer measure $1$ and
$\DD^{\TTT^\cccc}_\FF \cap P = \emptyset$
(since $\CC^\TTT_B \ne \TTT$ for all $B$),
so $\DD^{\TTT^\cccc}_\FF$ has inner measure $0$.
\end{proof}

\newpage
\section{Appendix}
\label{sec-app}

Here we collect a few proofs which don't seem worth putting
in the published part of the paper, since they just verify
some remarks in Section \ref{sec-ex}, which is itself essentially
an appendix.

First, as remarked in Section \ref{sec-ex}, one can get
a P-point in a modified version of Proposition \ref{prop-filter}.

\begin{proposition}
\label{prop-P-pt}
Assume Martin's Axiom (or, just $\pppp = \cccc$).
Suppose that $X$ is compact,  $Y$ is metric, and
$g$ is a limit point of $\langle f_n : n \in \omega\rangle$ in
$C_p(X,Y)$. 
Then there is a P-point ultrafilter
$\FF \subset [\omega]^\omega$ such that 
for each $x \in X$ there is 
$B \in \FF$ such that $\langle f_n(x) : n\in B\rangle
\to g(x)$.
\end{proposition}
\begin{proof}
As remarked in Section \ref{sec-ex}, it is sufficient to 
construct a P-point extending the
induced neighborhood filter, $\NN$.
As usual, $\FF$ will be generated by sets $\{B_\alpha : \alpha < \cccc\}$,
where $\alpha < \beta \to B_\beta \subseteq^* B_\alpha$.
Make sure, inductively, that
$$
g \mbox{ is a limit point of } \langle f_n : n \in B_\alpha \rangle
\ \ .  \eqno{(\mbox{\ding{166}})}
$$
To make $\FF$ an ultrafilter, list $[\omega]^\omega$
as $\{A_\alpha : \alpha < \cccc\}$, and 
let $B_{\alpha+1}$ be either
$B_\alpha \cap A_\alpha$ or $B_\alpha \backslash A_\alpha$,
preserving (\mbox{\ding{166}}), which ensures in particular that
if $A_\alpha \in \NN$ then $A_\alpha \in \FF$.

Now, say $\gamma$ is a limit.  We need $B_\gamma \subseteq^* B_\alpha$
for all $\alpha < \gamma$, along with (\ding{166}) for $B_\gamma$.
When $1 \le k < \omega$, let $\SS_k$ be the set of all finite
$S \subseteq \omega$ such that $\min(S) > k$ and
$$
\forall x_0, \ldots x_{k-1}\in X\;  \exists n \in S \;
\forall i < k\;  [ d(f_n(x_i), g(x_i)) < 1/k]\ \ .  \eqno{(\mbox{\ding{76}})}
$$
Observe that (\ding{166}) will hold if for all $k$, there
is an $S \subseteq B_\gamma$ with $S \in \SS_k$.
As in the proof of Proposition \ref{prop-filter}, there
is an
$S \subseteq B_{\alpha_1} \cap \cdots \cap  B_{\alpha_\ell}$
with $S \in \SS_k$ whenever $\alpha_1 < \cdots < \alpha_\ell < \gamma$.
Now, we obtain an appropriate $B_\gamma$
by a standard application of Martin's Axiom.
\end{proof}

Next, the argument in Proposition \ref{proposition-NN-rat} requires:

\begin{lemma}
Consider $\omega \subset \ZZZ \subseteq \bohr \ZZZ = \Hom(\TTT_d, \ZZZ)$,
and let $\lambda$ be Haar measure on $\bohr \ZZZ$. Fix $B \subseteq \omega$.
Then:
\begin{itemizz}
\item[1.]  If $B \in \NNb$, then $\lambda (\overline B) > 0$.
\item[2.]  If  $w \in \TTT$ has infinite order then
$ \langle w^n : n \in B\rangle \to 1$,
then $\lambda (\overline B) = 0$.
\end{itemizz}
\end{lemma}
\begin{proof}
For (1): We can assume that $B = U_1 \cap \cdots  \cap U_k \cap \omega$,
where $U_j = \{n \in \ZZZ : d(z_j^n, 1) < \varepsilon\}$.
Then $B \cup -B =  U_1 \cap \cdots \cap  U_k$, a basic open set in $\ZZZ^\#$.
So, $\lambda (\overline{B \cup -B}) > 0$
because $\ZZZ$ is dense in $\bohr \ZZZ$.
But $\overline{B \cup -B} = \overline B \cup \overline{-B}$ and
$\lambda(\overline B) = \lambda( \overline{-B})$, so 
$\lambda(\overline B) > 0$.

For (2):  Fix $\varepsilon > 0$; we show that
$\lambda (\overline B)  \le 2 \varepsilon$.
Let $d$ be the metric on $\TTT$ used in the proof of Lemma
\ref{lemma-limits-in-solenoids}.  By removing finitely many
elements of $B$, we can assume, WLOG, that 
$d(w^n, 1) \le \varepsilon$ for all $n\in B$.
Identifying $\bohr \ZZZ$ with $\Hom(\TTT_d, \ZZZ)$, let
$F = \{\varphi \in \bohr \ZZZ : d(\varphi(w), 1) \le \varepsilon\}$.
Applying duality (viewing $\TTT_d$ as $\widehat{\bohr\ZZZ}$),
let $E_w(\varphi) = \varphi(w)$; then $E_w$ is a continuous
homomorphism from $\bohr \ZZZ$ into $\TTT$, and
$F = \{\varphi \in \bohr \ZZZ : d(E_w(\varphi), 1) \le \varepsilon\}$,
which is a closed subset of $\bohr\ZZZ$.  Also, 
identifying $\omega \subset \ZZZ \subseteq \bohr \ZZZ$, we have
$B \subseteq F$, so it is sufficient to show that 
$\lambda (F)  = 2 \varepsilon$ if $\varepsilon \le 1/2$.
But this follows from the fact that $E_w$ maps \textit{onto} $\TTT$ (since
$w$ has infinite order), since $\{x \in \TTT : d(x,1) \le \varepsilon\}$
has measure $ 2 \varepsilon$.
\end{proof}

\end{document}